\documentclass{article}%
\usepackage{amssymb}
\usepackage{amsmath}
\usepackage{amsfonts}
\usepackage{graphicx}%
\providecommand{\U}[1]{\protect\rule{.1in}{.1in}}
\begin{document}

\begin{center}
\textbf{RULED SURFACES OF FINITE TYPE IN\ }

\textbf{3-DIMENSIONAL HEISENBERG\ GROUP}

\bigskip

Mohammed BEKKAR.
\end{center}

\bigskip

\textbf{Abstract. }In this paper, on the first, we prove $\Delta r=2$H where
$\Delta$ is the Laplacian operator, $r=\left(  r_{1},r_{2},r_{3}\right)  $ the
position vector field and H\textit{\ }is the mean curvature vector
field\textit{ }of a surface $\mathcal{S}$\ in the 3-dimensional Heisenberg
group $\mathbb{H}_{3}.$ In the second, we classify the ruled surfaces by
straight geodesic lines, which are of finite type in $\mathbb{H}_{3}.$ The
straight geodesic lines belong to $\ker\omega,$ where $\omega$ is the Darboux form.

\bigskip



\section{Introduction}

\subsection{Submanifolds of finite type in $3$-dimensional space}

Finite type submanifolds were introduced by B.-Y. Chen. A submanifold
$\mathcal{M}^{n}$ of an Euclidean space $\mathbb{E}^{n+p}$ is said to be of
finite type if each component of its position vector field $r$ can be written
as a finite sum of eigenfunctions of the Laplacian $\Delta$ of $\mathcal{M}%
^{n}$, i.e. if $r=r_{o}+r_{1}+r_{2}+...+r_{k}$ where $r_{o}$ is a constant and
$r_{1},r_{2},...,r_{k}$ non constant maps such that $\Delta r_{i}=\lambda
_{i}r_{i},\lambda_{i}\in\mathbb{R},i=1,..,k.$ If, all eigenvalues $\lambda
_{i},i=1,..,k$ are different, $\mathcal{M}^{n}$ is said to be of $k$-type.

The interesting examples of finite type submanifolds are some surfaces of
3-dimensional spaces of the Euclidean space $\mathbb{E}^{3}$.

During $1980-1990$, $\left[  6,7,8,9\right]  $ B.-Y. Chen introduced the
notion of finite type immersion in the dimensional Euclidean space
$\mathbb{E}^{m},$ in the pseudo Euclidean space and in the riemannian
manifolds. In the Euclidean space $\mathbb{E}^{3},$ this notion is, in one
way, a natural extension of the notion of minimal surfaces, which will be
itself, an extension of totaly geodesic surfaces\textit{.}

If we represent a such hypersurface in $\mathbb{E}^{m}$ for $m=3$\ we have in
$\mathbb{E}^{3}$

$\left(  1.1\right)  \qquad\qquad\qquad\mathcal{S}:r(x,y)=\left(
r_{1}(x,y),r_{2}(x,y),r_{3}(x,y)\right)  ,(x,y)\in\mathcal{D}\subseteq
\mathbb{R}^{2},$

It is well known that in $\mathbb{E}^{3},$ for all regular surface
$\mathcal{S}$

$\left(  1.2\right)  \qquad\qquad\qquad\qquad\qquad\qquad\Delta r=-2$H

where $\Delta$ is the Laplace operator and H\textit{\ }is the mean curvature
vector field\textit{ }of $\mathcal{S}$ The components $\left(  r_{1}%
,r_{2},r_{3}\right)  $ of $r$ are differentiable functions so that

$\left(  1.3\right)  \qquad\qquad\qquad\qquad\Delta r(x,y)=\left(  \Delta
r_{1}(x,y),\Delta r_{2}(x,y),\Delta r_{3}(x,y)\right)  $

From $\left(  1.2\right)  ,$\ minimal surfaces and spheres verify $\Delta
$H$=\lambda$H$,\lambda\in\mathbb{R}.$ $\left(  1.2\right)  $\ shows that
$\mathcal{S}$\ is minimal surface of $\mathbb{E}^{3},$\ if and only if,
$r_{i},i=1,2,3$ are harmonic.

If the matrix $\left(  g_{ij}\right)  $ consists on the components of the
induced any metric on $\mathcal{S}$, and $\left(  g^{ij}\right)  $ its inverse
and $D=\det\left(  g_{ij}\right)  $; the Laplacian (Beltrami's operator)
$\Delta$\ on $\mathcal{S}$ is given by

$\left(  1.4\right)  $\qquad$\qquad\qquad\qquad\Delta=\frac{-1}{\sqrt
{\left\vert D\right\vert }}%
%TCIMACRO{\dsum \limits_{i,j=1}^{n}}%
%BeginExpansion
{\displaystyle\sum\limits_{i,j=1}^{n}}
%EndExpansion
\frac{\partial}{\partial x^{i}}\left(  \sqrt{\left\vert D\right\vert }%
g^{ij}\frac{\partial}{\partial x^{j}}\right)  .$

The classifcation of $1$-type submanifolds of Euclidean space $\mathbb{E}^{3}$
was done in 1966 by T. Takahashi $[21]$. He proved that the submanifolds in
$\mathbb{R}^{m}$ satisfy the differential equation

$(1.5)$\qquad$\qquad\qquad\qquad\qquad\Delta r=\lambda r,$

for some real number $\lambda$, if and only if, either the submanifold is a
minimal submanifold of $\mathbb{R}^{m},(\lambda=0),$ or a hypersphere of
$\mathbb{R}^{m}$ centered at the origin $(\lambda\neq0)$. O. J. Garay $[15]$
generalized $\left(  1.5\right)  $ where he studied hypersurfaces in
$\mathbb{R}^{m}$, not necessarily associated to the same eigenvalue. He linked
and considered hypersurfaces in $\mathbb{R}^{m}$ satisfying the differential equation

$(1.6)\qquad\qquad\qquad\qquad\Delta r=Ar,A\in M(m,\mathbb{R})$.

Next, F. Dillen, J. Pas and L. Verstraelen $\left[  10\right]  $ observed and
proposed the study of submanifolds of $\mathbb{R}^{m}$ satisfying the
following equation

\begin{center}
$\Delta r=Ar+B,B\in\mathbb{R}^{m}.$
\end{center}

The same and others authors studied several problems on link to the subject of
finite type particular surfaces like translation surfaces, the quadrics,
surfaces of revolution, helicoidal surfaces,...

The study of this notion of finite type was extended for surfaces in Euclidean
$3$-space $\mathbb{E}^{3}$\ of specific form such that their Gauss map
$\boldsymbol{N}$ satisfies an analogous or similar equation $\Delta
\boldsymbol{N}=A\boldsymbol{N}\mathbf{,}A\in M(m,\mathbb{R})$.

In the same way, many authors studied particular surfaces of finite type in
the Euclidean, pseudo Euclidean and Lorentz-Minkowski 3-dimensional space
satisfying the differential equation $\Delta^{II}r=Ar;\Delta^{III}r=Ar$, where
$\Delta^{II},\ \Delta^{III}$ are respectively the Laplace operator with
respect to the second and third fundamental form which are not degenerated.
For these works of different authors, see the references and therein. Among
from others, we have B-.Y. Chen, L. J. Alias, A. Ferrandez, C. Baikoussis, D.
E. Blair, Piccinni, S. M. Choi, Kim Y. H., Yoon, D. W, P. Lucas... The list of
authors working in the subject is very long and certainly not closed.

In our work, first, we proof that $\left(  1.2\right)  $\ stays true in
Heisenberg space $\mathbb{H}_{3}$ and classify ruled surfaces by the straight
lines as geodesics of the Heisenberg group $\mathbb{H}_{3}$ which are of
finite type . Explicitly, we search ruled surfaces by the straight lines
geodesics of $\mathbb{H}_{3}$\ which satisfy

$(1.7)$\qquad$\qquad\qquad\qquad\qquad\Delta r_{i}=\lambda_{i}r_{i}%
,\lambda_{i}\in\mathbb{R},i=1,2,3.$

\subsection{Minimal surfaces in 3-dimensional space.}

Let $z=f(x,y)$ be a graph of a regular surface $\mathcal{S}$ in the Euclidean
space $\mathbb{E}^{3}.$ $\mathcal{S}$ is minimal surface if $f$ satisfies

$(1.8)$\qquad$\qquad f_{xx}\left(  1+f_{y}^{2}\right)  -2f_{x}f_{y}%
f_{xy}+f_{yy}\left(  1+f_{x}^{2}\right)  =0,f_{x}=\frac{\partial f}{\partial
x}....$

This equation obtained by J. L. Lagrange in 1776 where he used the variational
calculus. In 1775, J. B. Meusnier gave a geometric interpretation to this,
said that $\mathcal{S}$\ is minimal if its mean curvature function $H\equiv0.$
Some particular solutions have been given by the same or others authors until
1866 when Weirstrass\textit{\ }solved $(1.8).$

The author construct examples of minimal surfaces with the choise of two
complex functions so that the parameters in representation of $\mathcal{S}%
$\ are isothermal.

Since, more mathematicians work and study different directions on minimal
surfaces and several questions therein, like to found minimal surfaces in
other spaces, like in Lorentz space raised but unsuccessfuly..., until in 1982
when Do Carmo M. found an original and beautiful minimal surface.

S. N. Bernstein, proved $\left(  1915-1917\right)  $ that \textit{" there are
no other complete graphs, except the planes, which are minimal surfaces in
}$\mathbb{E}^{3}$.

In 1990, a natural question raised for the $\mathbb{H}_{3},$ are there, as
$\mathbb{E}^{3}$, minimal surfaces in $\mathbb{H}_{3}?$

This space is the 3-dimensional Heisenberg group which can be seen as
$\mathbb{R}^{3}$ equiped with a riemannian metric

\begin{center}
$ds_{\mathbb{H}_{3}}^{2}=dx^{2}+dy^{2}+\left(  dz+\frac{1}{2}(ydx-xdy)\right)
^{2}.$
\end{center}

To know deeply the geometry of $\mathbb{H}_{3}$, in 1991, the author $\left[
2\right]  $ searched and wrote the equation of minimal surfaces as a graph of
functions $z=f(x,y)$

\begin{center}
$(1.9)$

$f_{xx}\left(  1+\left(  f_{y}-\frac{x}{2}\right)  ^{2}\right)  -2f_{xy}%
\left(  f_{x}+\frac{y}{2}\right)  \left(  f_{y}-\frac{x}{2}\right)
+f_{yy}\left(  1+\left(  f_{x}+\frac{y}{2}\right)  ^{2}\right)  =0,$
\end{center}

and gave some particular solutions $\left[  2\right]  $. Among these, ones are
already minimal surfaces in $\mathbb{E}^{3}$ like planes, helico\"{\i}ds. The
other particular solution $f(x,y)=\frac{xy}{2}$\ of $(1.9),$ which is not
solution of $(1.8),$ is the \textit{hyperbolic paraboloid }which have a very
important part as a minimal surface in $\mathbb{H}_{3}$.

In an other work, T. Sari together with the author $\left[  3\right]  $ gave a
complete description of minimal surfaces ruled by straight lines as geodesics
of $\mathbb{H}_{3}$\ and lines. By the analogous method used by Weirstrass
process to solve $(1.8)$, F. Mercuri, S. Montaldo and P. Piu solved $(1.9),$
$\left[  20\right]  $. Other works of different authors followed. Bernstein's
theorem is not valid in $\mathbb{H}_{3}$\ since $z=\frac{xy}{2}$ is a complete
minimal surface, other the plane.

\section{Preleminaries}

\subsection{Heisenberg space}

$\bullet$ Heisenberg group $\mathbb{H}_{3}$\ is a two-step nilpotent Lie group
which is a subgroup of linear group $Gl(3,\mathbb{R}).$ It is known as a
quantic physic model strongly studied by the theoritical physicists and
mathematicians. It is a 3-dimensional riemannian manifold equiped with the
riemannian metric

$(2.1)$\qquad$\qquad\qquad\qquad\qquad ds_{\mathbb{H}_{3}}^{2}=dx^{2}%
+dy^{2}+\omega^{2}$

where $\omega=dz+\frac{1}{2}(ydx-xdy),$ is an Pfaffian form, known as Darboux
form, it is also an contact form.

$\bullet$ As a Lie group, $\mathbb{H}_{3}$\ acts by left translations keeping
invariant $ds_{\mathbb{H}_{3}}^{2}$

\begin{center}
$\mathbb{H}_{3}=\left\{  \left(
\begin{array}
[c]{ccc}%
1 & x & z\\
0 & 1 & y\\
0 & 0 & 1
\end{array}
\right)  ,\left(  x,y,z\right)  \in\mathbb{R}^{3}\right\}  .$
\end{center}

The one to one map

\begin{center}
$\mu:\mathbb{R}^{3}\longrightarrow\mathbb{H}_{3}$

$\left(
\begin{array}
[c]{c}%
x\\
y\\
z
\end{array}
\right)  \longrightarrow\left(
\begin{array}
[c]{ccc}%
1 & x & z\\
0 & 1 & y\\
0 & 0 & 1
\end{array}
\right)  $
\end{center}

induces a group structure.

$\bullet$ $\mathbb{H}_{3}$ is equipped with the no commutatif group structure
$\left(  \ast\right)  $\ given by

$(2.2)$\qquad$\qquad\left(  x,y,z\right)  \ast\left(  x^{\prime},y^{\prime
},z^{\prime}\right)  =\left(  x+x^{\prime},y+y^{\prime},z+z^{\prime}+\frac
{1}{2}(x^{\prime}y-xy^{\prime})\right)  ,$

\begin{center}
$\forall\left(  x,y,z\right)  ,\left(  x^{\prime},y^{\prime},z^{\prime
}\right)  \in\mathbb{R}^{3}.$
\end{center}

$\bullet$ These 1-forms $dx,dy,\omega$ in $ds_{\mathbb{H}_{3}}^{2}$\ are
invariant by left translations in $\mathbb{H}_{3}$ and by rotations about
$(Oz)$ axis. The left invariant orthonormal coframe is associate with the
orthonormal left invariant frame

\begin{center}
$e_{1}=\partial_{x}-\frac{y}{2}\partial_{z},e_{2}=\partial_{y}+\frac{x}%
{2}\partial_{z},e_{3}=\partial_{z}.,\partial_{x}=\frac{\partial}{\partial x}.$
\end{center}

$\left\{  \partial_{x},\partial_{y},\partial_{z}\right\}  $ denote the natural
vector fields in $\mathbb{R}^{3}.$ The Lie brackets are

\begin{center}
$\left[  e_{2},e_{3}\right]  =\left[  e_{3},e_{1}\right]  =0,\left[
e_{1},e_{2}\right]  =e_{3}.$
\end{center}

$\bullet$ The Levi-Civita connection denoted $\tilde{\nabla}$ of
$\mathbb{H}_{3}$ is given by

\begin{center}
$\left(
\begin{array}
[c]{c}%
\tilde{\nabla}_{e_{1}}e_{1}\\
\tilde{\nabla}_{e_{1}}e_{2}\\
\tilde{\nabla}_{e_{1}}e_{3}%
\end{array}
\right)  =\left(
\begin{array}
[c]{c}%
0\\
\frac{1}{2}e_{3}\\
-\frac{1}{2}e_{2}%
\end{array}
\right)  ,\left(
\begin{array}
[c]{c}%
\tilde{\nabla}_{e_{2}}e_{1}\\
\tilde{\nabla}_{e_{2}}e_{2}\\
\tilde{\nabla}_{e_{2}}e_{3}%
\end{array}
\right)  =\left(
\begin{array}
[c]{c}%
-\frac{1}{2}e_{3}\\
0\\
\frac{1}{2}e_{1}%
\end{array}
\right)  ,$

$\left(
\begin{array}
[c]{c}%
\tilde{\nabla}_{e_{3}}e_{1}\\
\tilde{\nabla}_{e_{3}}e_{2}\\
\tilde{\nabla}_{e_{3}}e_{3}%
\end{array}
\right)  =\left(
\begin{array}
[c]{c}%
-\frac{1}{2}e_{2}\\
\frac{1}{2}e_{1}\\
0
\end{array}
\right)  .$
\end{center}

$\bullet$ Let $\mathcal{S}$ in $\mathbb{H}_{3}$ be an orientable surface and
$r:\mathcal{S}\longrightarrow\mathbb{H}_{3}$ an immersion. Denote by $\nabla$
the induced Levi-Civita connection on $\mathcal{S}$. If $\boldsymbol{N}$ is
the unit normal vector on $\mathcal{S}$, we have the well known Gauss and
Weingarteen formulae for riemannian manifold and hypersurfaces

\begin{center}
$\left\{
\begin{array}
[c]{c}%
\tilde{\nabla}_{X}Y=\nabla_{X}Y+h(X,Y)\boldsymbol{N}\mathbf{,}h(X,Y)=g\left(
\tilde{\nabla}_{X}Y,\boldsymbol{N}\right) \\
\tilde{\nabla}_{X}\boldsymbol{N}=-\mathcal{A}X\text{
\ \ \ \ \ \ \ \ \ \ \ \ \ \ \ \ \ \ \ \ \ \ \ \ \ \ \ \ \ \ \ \ \ \ \ \ \ \ \ \ \ \ \ \ \ \ \ }%
\end{array}
\right.  $
\end{center}

where $X,Y$ are tangent vectors fields on $\mathcal{S}$, $h$ the second
fundamental form and $A$ is the shape operator associated to $\boldsymbol{N}$.
Recall that $\mathcal{A}$ is a self adjoint endomorphism with respect to the
metric on $\mathcal{S}$, that is

\begin{center}
$g(\mathcal{A}X,Y)=g(X,\mathcal{A}Y),X,Y\in T\mathcal{S}$.
\end{center}

The mean curvature of the immersion $r$ is H$=\frac{1}{2}tr(\mathcal{A}).$

$\bullet$ $\mathbb{H}_{3}$ is a homogenous Riemannian manifold after Euclidean
space $\mathbb{E}^{3},$\ the sphere $\mathbb{S}^{3}$ and hyperbolic space
$\mathbb{H}^{3}$. Its isometry group is 4-dimensional and it is isomorphic to
semi direct sum $SO(2)\vartriangleright\mathbb{R}^{3}.$ Explicitly, the
component of identity is the affine group of transformations of $\mathbb{R}%
^{3}$

\begin{center}
$\left(
\begin{array}
[c]{c}%
x\\
y\\
z
\end{array}
\right)  \longrightarrow\left(
\begin{array}
[c]{ccc}%
\cos\theta & -\sin\theta & 0\\
\sin\theta & \cos\theta & 0\\
A & B & 1
\end{array}
\right)  \left(
\begin{array}
[c]{c}%
x\\
y\\
z
\end{array}
\right)  +\left(
\begin{array}
[c]{c}%
a\\
b\\
c
\end{array}
\right)  $
\end{center}

where $A=\frac{1}{2}\left(  a\sin\theta-b\cos\theta\right)  ,B=\frac{1}%
{2}\left(  a\cos\theta+b\sin\theta\right)  .$

The isometry group of $\mathbb{H}_{3}$ is discribed by $\left\{  \left(
\theta,a,b,c\right)  ;\theta,a,b,c\in\mathbb{R}\right\}  $. It contains the
rotations $\left(  \theta,0,0,0\right)  $\ about $(Oz)$ and the left
translations $\left(  0,a,b,c\right)  .$

\textbf{Remarks}

1. The left translations of the isometry group of $\mathbb{H}_{3}$\ don't keep
left invariant Euclidean metrics$.$\ As the same, the Euclidean translations
in $\mathbb{H}_{3}$\ don't keep invariant Heisenberg metrics $ds_{\mathbb{H}%
_{3}}^{2}$.

2. The determinant of matrix of $\mathbb{H}_{3},$ and the one associated to
the metric is constant equal to one. It is an important characteristic of
$\mathbb{H}_{3}$\ and especially for its Riemannian metric $ds_{\mathbb{H}%
_{3}}^{2}.$

\bigskip

$\bullet$ As in $\mathbb{E}^{3},$ let $z=f(x,y)$ be a graph of a regular
surface $\mathcal{S}$ in $\mathbb{H}_{3}$ with its position vector field

\begin{center}
$r(x,y)=(x,y,f(x,y));(x,y)\in\mathcal{D}\subseteq\mathbb{R}^{2}.$
\end{center}

The first fundamental form of $\mathcal{S}$ is obtained like trace of
$ds_{\mathbb{H}_{3}}^{2}$ on $\mathcal{S}$,

\begin{center}
$ds_{\mathbb{H}_{3\mid\mathcal{S}}}^{2}=g_{\mid_{\mathcal{S}}}=dx^{2}%
+dy^{2}+\left(  f_{x}dx+f_{y}dy+\frac{1}{2}\left(  ydx-xdy\right)  \right)
^{2}$

$=\left(  1+P^{2}\right)  dx^{2}+2PQdxdy+\left(  1+Q^{2}\right)  dy^{2}.$
\end{center}

where $P=f_{x}+\frac{y}{2},Q=f_{y}-\frac{x}{2}.$

$\bullet$ A basis of tangent space $\mathcal{T}_{p}\mathcal{S}$ on
$p\in\mathcal{S}$ associated to

\begin{center}
$\mathcal{S}:r\left(  x,y\right)  =\left(  x,y,f\left(  x,y\right)  \right)  $
\end{center}

is given by

$\left(  2.3\right)  \qquad\qquad\qquad\qquad\left\{
\begin{array}
[c]{c}%
r_{x}=(1,0,f_{x})=\partial_{x}+f_{x}\partial_{z}=e_{1}+Pe_{3}\\
r_{y}=(0,1,f_{y})=\partial_{y}+f_{y}\partial_{z}=e_{2}+Qe_{3}%
\end{array}
\right.  $

The coefficients of the first fundamental form of $\mathcal{S}$ are

$(2.4)\qquad\left\{
\begin{array}
[c]{c}%
E=g(r_{x},r_{x})=1+P^{2},G=g(r_{y},r_{y})=1+Q^{2}\text{ \ \ \ \ \ \ }\\
F=g\left(  r_{x},r_{y}\right)  =PQ,;EG-F^{2}=1+P^{2}+Q^{2}=W^{2}%
\end{array}
\right.  $

$\bullet$ The unit normal vector field $\boldsymbol{N}$\ on\ $\mathcal{S}$ is

\begin{center}
$\boldsymbol{N}=\frac{1}{W}\left(  Pe_{1}+Qe_{2}-e_{3}\right)  .$
\end{center}

In order to compute the coefficients of the second fundamental form of
$\mathcal{S}$, we have to calculate the following

$\left(  2.5\right)  $\qquad$\qquad\ \ \ \left\{
\begin{array}
[c]{c}%
r_{xx}=\tilde{\nabla}_{r_{x}}r_{x}=-Pe_{2}+f_{xx}e_{3}\text{
\ \ \ \ \ \ \ \ \ \ \ \ \ \ \ \ \ \ }\\
r_{xy}=\tilde{\nabla}_{r_{x}}r_{y}=\tilde{\nabla}_{r_{y}}r_{x}=\frac{P}%
{2}e_{1}-\frac{Q}{2}e_{2}+f_{xy}e_{3}\\
r_{yy}=\tilde{\nabla}_{r_{y}}r_{y}=Qe_{1}+f_{yy}e_{3}\text{
\ \ \ \ \ \ \ \ \ \ \ \ \ \ \ \ \ \ \ \ \ }%
\end{array}
\right.  $

wich imply that the coefficients of the second fundamental form of
$\mathcal{S}$ are given by

$(2.6)$\qquad$\qquad\left\{
\begin{array}
[c]{c}%
L=-g(\tilde{\nabla}_{r_{x}}r_{x},\boldsymbol{N})=\frac{f_{xx}+PQ}{W}\text{
\ \ \ \ \ \ \ \ }\\
M=-g(\tilde{\nabla}_{r_{x}}r_{y},\boldsymbol{N})=\frac{f_{xy}+\frac{1}{2}%
Q^{2}-\frac{1}{2}P^{2}}{W}\\
N=-g(\tilde{\nabla}_{r_{y}}r_{y},\boldsymbol{N})=\frac{f_{yy}-PQ}{W}\text{
\ \ \ \ \ \ \ \ }%
\end{array}
\right.  .$

In $\left[  3\right]  $ we have some expressions which simplify the computation

$(2.7)\qquad\left\{
\begin{array}
[c]{c}%
P_{x}=f_{xx},Q_{y}=f_{yy},P_{y}=f_{xy}+\frac{1}{2},Q_{x}=f_{xy}-\frac{1}%
{2},P_{y}-Q_{x}=1\\
P_{y}+Q_{x}=2f_{xy};Q=-uP,u^{2}P_{x}+uP_{y}+uQ_{x}+Q_{y}=0\text{ \ \ \ \ \ }%
\end{array}
\right.  $

If we put $H_{1}=EN-2FM+GL$, the mean curvature $H$ of $\mathcal{S}$, with the
help of $\left(  2.7\right)  ,$ is

$(2.8)\qquad\qquad\qquad H=\frac{H_{1}}{2W^{2}},H_{1}=\frac{1}{W}\left(
P_{x}+Q_{y}\right)  =\frac{1}{W}\left(  f_{xx}+f_{yy}\right)  .$

\subsection{Characterization of minimal surfaces in $\mathbb{H}_{3}$}

We proved in $\left[  3\right]  ,$\ the theorem 2$.$ \textit{The minimal
surfaces of }$\mathbb{H}_{3}$\ \textit{ruled by straight geodesic lines are
localy the graph of harmonic functions.} In the next theorem, we prove,
$\mathcal{S}$\textit{\ }is a minimal surfaces in $\mathbb{H}_{3},$ if and only
if its coordinate functions are harmonic.

\subsubsection{Theorem}

\textit{The Beltrami formula in }$\mathbb{H}_{3}$\textit{\ is given by
}$\Delta r=2$H$,$\textit{\ }$\Delta$\textit{\ is the Laplacian of the surface
with respect the first fundamental form of the surface }$\mathcal{S}%
:r(x,y)$\textit{\ in }$\mathbb{H}_{3},$\textit{\ }H\textit{\ is the mean
curvature vector field of }$\mathcal{S}.$

\bigskip

\textbf{Proof}

We have explicitly the Beltrami operator $\left(  1.4\right)  $

$(2.9)\qquad\Delta=\frac{-1}{W}\left(  \frac{\partial}{\partial x}\left(
Wg^{11}\frac{\partial}{\partial x}+Wg^{12}\frac{\partial}{\partial y}\right)
+\frac{\partial}{\partial y}\left(  Wg^{21}\frac{\partial}{\partial x}%
+Wg^{22}\frac{\partial}{\partial y}\right)  \right)  .$

The matrix associated to the first fundamental form of $\mathcal{S}$\ is

\begin{center}
$\left(  g_{ij}\right)  =\left(
\begin{array}
[c]{cc}%
E & F\\
F & G
\end{array}
\right)  ,\left(  g^{ij}\right)  =\left(
\begin{array}
[c]{cc}%
\frac{G}{EG-F^{2}} & -\frac{F}{EG-F^{2}}\\
-\frac{F}{EG-F^{2}} & \frac{E}{EG-F^{2}}%
\end{array}
\right)  ,\det\left(  g_{ij}\right)  =W^{2}.$
\end{center}

The operator $\Delta$\ turns out to

\begin{center}
$\Delta=\frac{-1}{W}\left[  \left(  \frac{\partial}{\partial x}\left(
\frac{G}{W}\frac{\partial}{\partial x}-\frac{F}{W}\frac{\partial}{\partial
y}\right)  \right)  +\left(  \frac{\partial}{\partial y}\left(  -\frac{F}%
{W}\frac{\partial}{\partial x}+\frac{E}{W}\frac{\partial}{\partial y}\right)
\right)  \right]  .$
\end{center}

The first part of $\Delta$

\begin{center}
$\Delta_{1}=\frac{-1}{W}\left(  \frac{\partial}{\partial x}\left(  \frac{G}%
{W}\frac{\partial}{\partial x}-\frac{F}{W}\frac{\partial}{\partial y}\right)
\right)  $

$=\frac{-1}{W}\left(
\begin{array}
[c]{c}%
\frac{2QQ_{x}W-G\left(  PP_{x}+QQ_{x}\right)  W^{-1}}{W^{2}}\frac{\partial
}{\partial x}-\frac{\left(  P_{x}Q+PQ_{x}\right)  W-F\left(  PP_{x}%
+QQ_{x}\right)  W^{-1}}{W^{2}}\frac{\partial}{\partial y}\\
+\frac{G}{W}\frac{\partial^{2}}{\partial x^{2}}-\frac{F}{W}\frac{\partial^{2}%
}{\partial x\partial y}%
\end{array}
\right)  .$
\end{center}

The second part

\begin{center}
$\Delta_{2}=\frac{-1}{W}\left(  \frac{\partial}{\partial y}\left(  -\frac
{F}{W}\frac{\partial}{\partial x}+\frac{E}{W}\frac{\partial}{\partial
y}\right)  \right)  $

$=\frac{-1}{W}\left(
\begin{array}
[c]{c}%
-\frac{\left(  P_{y}Q+PQ_{y}\right)  W-F\left(  PP_{y}+QQ_{y}\right)  W^{-1}%
}{W^{2}}\frac{\partial}{\partial x}+\frac{2PP_{y}W-E\left(  PP_{y}%
+QQ_{y}\right)  W^{-1}}{W^{2}}\frac{\partial}{\partial y}\\
+\frac{E}{W}\frac{\partial^{2}}{\partial y^{2}}-\frac{F}{W}\frac{\partial^{2}%
}{\partial x\partial y}%
\end{array}
\right)  .$
\end{center}

Thus, $\Delta$\ turns out to

\begin{center}
$\Delta=-\frac{1}{W^{2}}\left(
\begin{array}
[c]{c}%
G\frac{\partial^{2}}{\partial x^{2}}-2F\frac{\partial^{2}}{\partial x\partial
y}+E\frac{\partial^{2}}{\partial y^{2}}\\
-\left(  Q+PQ_{y}-QQ_{x}\right)  \frac{\partial}{\partial x}+\left(
-P+QP_{x}-PP_{y}\right)  \frac{\partial}{\partial y}%
\end{array}
\right)  $

$+\frac{1}{W^{4}}\left(
\begin{array}
[c]{c}%
\left(  \left(  PP_{x}+QQ_{x}\right)  G-\left(  PP_{y}+QQ_{y}\right)
F\right)  \frac{\partial}{\partial x}\\
+\left(  \left(  PP_{y}+QQ_{y}\right)  E-\left(  PP_{x}+QQ_{x}\right)
F\right)  \frac{\partial}{\partial y}%
\end{array}
\right)  .$
\end{center}

With the help of $\left(  2.3\right)  -\left(  2.7\right)  $

\begin{center}
$\Delta=-\frac{1}{W^{4}}\left(
\begin{array}
[c]{c}%
\left(  FPP_{y}-GPP_{x}+FQQ_{y}-GQQ_{x}\right)  e_{1}\\
+\left(  FPP_{x}+FQQ_{x}-EPP_{y}-EQQ_{y}\right)  e_{2}\\
+\left(
\begin{array}
[c]{c}%
FP^{2}P_{y}-GP^{2}P_{x}+FQ^{2}Q_{x}-Q^{2}EQ_{y}\\
+FPQP_{x}+FPQQ_{y}-GPQQ_{x}-PQEP_{y}%
\end{array}
\right)  e_{3}%
\end{array}
\right)  .$
\end{center}

The computation gives

\begin{center}
$(2.10)$\qquad$\Delta r=2H\boldsymbol{N=}\frac{1}{W^{3}}\left(  P_{x}%
+Q_{y}+Q^{2}P_{x}+P^{2}Q_{y}-PQ-2PQQ_{x}\right)  =2$H$.$
\end{center}

\bigskip

\textbf{Remark}

We can, also write the Beltrami operator $\Delta$ in the form

\begin{center}
$\Delta=\frac{-1}{W^{4}}\left(
\begin{array}
[c]{c}%
W^{2}\left(  G\frac{\partial^{2}}{\partial x^{2}}-2F\frac{\partial^{2}%
}{\partial x\partial y}+E\frac{\partial^{2}}{\partial y^{2}}\right) \\
-\left(  PH_{1}+QW^{2}\right)  \frac{\partial}{\partial x}+\left(
-QH_{1}+PW^{2}\right)  \frac{\partial}{\partial y}%
\end{array}
\right)  .$
\end{center}

Thus, by a straightforward computation, and using $\left(  2.3\right)
-(2.8),$ we obtain

\begin{center}
$\Delta r=\frac{-1}{W^{2}}\left(  Gr_{xx}-2Fr_{xy}+Er_{yy}-Qr_{x}%
+Pr_{y}\right)  +\frac{1}{W^{4}}H_{1}\left(  Pr_{x}+Qr_{y}\right)  $
\end{center}

which give $(2.10).$

Hence $\left(  1.2\right)  $ is verified in $\mathbb{H}_{3}$. $\mathcal{S}%
$\textit{\ }is a minimal surface in $\mathbb{H}_{3},$ if only if the functions
$r_{i},i=1,2,3$ are harmonic. We conclude, as in the Euclidean space
$\mathbb{E}^{3}$, all the minimal surfaces in $\mathbb{H}_{3}$\ are of finite type.

\subsection{Ruled surfaces by geodesic lines in\ $\mathbb{H}_{3}.$}

In $\left[  17\right]  ,$\ Th. Hangan solved the Euler-Lagrange system of
geodesics in $\mathbb{H}_{3}$ and gave explicitly the solutions. Among these,
the author obtained a straight lines which are in $\ker\omega,$ $\omega$ is
the Darboux form.

After writting $\left(  1.9\right)  \ $and\ giving some particular\ solutions
in $\left[  2\right]  $, in 1992, T. Sari together with the author gave a
complete description of minimal surfaces ruled by straight lines as geodesics
and straight lines in $\mathbb{H}_{3}$ $\left[  3\right]  .$

In fact and to summarise, if we get the geodesic $\Gamma$ coming from a point
$p=(x(0),y(0),z(0))$ of $\mathbb{H}_{3}$\ and tangent to the vector
$v=(x^{\prime}(0),y^{\prime}(0),z^{\prime}(0))$ of $\mathcal{T}_{P}%
\mathbb{H}_{3},$ $\Gamma$ is an straight line in the two cases, when $v$ is
perpendicular which says that $\Gamma$\ is a perpendicular line or $v$ is in
$\ker\omega,$ saying that $\Gamma$ is a line in $\ker\omega$.

Explicitly, let a surface $\mathcal{S}$ in $\mathbb{H}_{3}$ parametrized by

\begin{center}
$(t,s)\longrightarrow\alpha(t)+s\beta(t),\alpha(t)\in\mathbb{R}^{3}%
,\beta(t)\in\mathbb{R}^{3}\smallsetminus\left\{  0\right\}  .$
\end{center}

$\mathcal{S}$ is ruled by the straight line $\mathcal{L}$ stem from the point
$\alpha(t)$ and has $\beta(t)$ as director vector.

Up to isometry of $\mathbb{H}_{3}$ and an appropriate choice of the curve
$\alpha(t),$ we obtain two families of ruled minimal surfaces by geodesic
lines $\left[  3\right]  $ which are

$\qquad\qquad\qquad\qquad\qquad\qquad1^{\circ})$ $\alpha(t)=(t,a(t),0),\beta
(t)=(0,0,1)$

$\qquad\qquad\qquad\qquad\qquad\qquad2^{\circ})$ $\alpha(t)=(t,0,a(t)),\beta
(t)=(u(t),1,\frac{t}{2})$

where $a(t)$ and $u(t)$ are differentiable functions of $t$.

For convenience, in the following, we used the notations of $\left[  3\right]
.$

\section{Ruled surfaces by geodesic lines satisfying $\Delta r_{i}=\lambda
_{i}r_{i}$\ in $\mathbb{H}_{3}.$}

\subsection{The first family of ruled surfaces by geodesic lines satisfying
$\Delta r_{i}=\lambda_{i}r_{i}$}

Take a regular surface $\mathcal{S}:r(t,s)$ in $\mathbb{H}_{3}$. The Laplacien
operator with respect to the first fundamental form for $\mathcal{S}$ is, for
each coordinate $r_{i},i=1,2,3$

$(3.1)$\qquad$\Delta r_{i}=-\frac{1}{W}\left(  \left(  \frac{Gr_{i,t}%
-Fr_{i,s}}{W}\right)  _{t}+\left(  \frac{-Fr_{i,t}+Er_{i,s}}{W}\right)
_{s}\right)  =\lambda_{i}r_{i}$

If we denote by $\mathcal{S}_{1}$ the first ruled surface

\begin{center}
$\mathcal{S}_{1}:r(t,s)=\left(  r_{1}(t,s),r_{2}(t,s),r_{3}(t,s)\right)
=(t,a(t),0)+s(0,0,1)=(t,a(t),s).$
\end{center}

we have $\left[  3\right]  $, by $(2.3),(2.6),$ the coefficients of the
fundamental forms

\begin{center}
$\left\{
\begin{array}
[c]{c}%
E=1+(a^{\prime})^{2}+\frac{1}{4}\left(  a-ta^{\prime}\right)  ^{2}%
,G=1,F=\frac{1}{2}\left(  a-ta^{\prime}\right)  ,W^{2}=1+(a^{\prime})^{2}\\
L=a^{\prime\prime}-\frac{1}{2}\left(  a-ta^{\prime}\right)  \left(
1+(a^{\prime})^{2}\right)  ,M=-\frac{1}{2}\left(  1+(a^{\prime})^{2}\right)
,N=0\text{. \ \ \ \ \ \ \ \ \ \ }%
\end{array}
\right.  $
\end{center}

The mean curvature of $\mathcal{S}_{1}:r(t,s)$\ is

\begin{center}
$H=\frac{EN+GL-2FM}{2(EG-F^{2})}=\frac{a^{\prime\prime}}{2\left(
1+(a^{\prime})^{2}\right)  }=\frac{a^{\prime\prime}}{2W^{2}}.$
\end{center}

The system $\left(  3.1\right)  $\ turns out for each $i=1,2,3$ to

$\bullet$ For $r_{1}\left(  t,s\right)  =t,$ $r_{1,t}=1,r_{1,t}=0,r_{1,s}%
=0,r_{1,s}=0,$

\begin{center}
$Gr_{1,t}-Fr_{1,s}=1,-Fr_{1,t}+Er_{1,s}=-\frac{1}{2}\left(  a-ta^{\prime
}\right)  ,$

$\Delta r_{1}(t,s)=\frac{2a^{\prime}a^{\prime\prime}}{2\left(  1+(a^{\prime
})^{2}\right)  ^{2}}=\frac{2a^{\prime}H}{W^{2}}.$
\end{center}

$\bullet$\textit{\ }For $r_{2}(t,s)=a(t)$, $r_{2,t}=a^{\prime};r_{2,s}=0,$

\begin{center}
$Gr_{2,t}-Fr_{2,s}=a^{\prime},-Fr_{2,t}+Er_{2,s}=-\frac{1}{2}\left(
a-ta^{\prime}\right)  a^{\prime},$

$\Delta r_{2}(t,s)=\frac{-a^{\prime\prime}}{\left(  1+(a^{\prime})^{2}\right)
^{2}}=-\frac{2H}{W^{2}}.$
\end{center}

$\bullet$\textit{\ }For $r_{3}(t,s)=s,$ $r_{3,t}=0,r_{3,s}=1,$

\begin{center}
$Gr_{3,t}-Fr_{3,s}=-\frac{1}{2}\left(  a-ta^{\prime}\right)  ,-Fr_{3,t}%
+Er_{3,s}=1+(a^{\prime})^{2}+\frac{1}{4}\left(  a-ta^{\prime}\right)  ^{2},$

$\Delta r_{3}(t,s)=\frac{-1}{2\left(  1+(a^{\prime})^{2}\right)  ^{2}%
}a^{\prime\prime}\left(  t+aa^{\prime}\right)  =\frac{-\left(  t+aa^{\prime
}\right)  H}{W^{2}}.$
\end{center}

Then the system $(1.7)$ with help of $(3.1)$ turns out to

$(3.2)$\qquad$\qquad\qquad\left\{
\begin{array}
[c]{c}%
\frac{2a^{\prime}H}{W^{2}}=\lambda_{1}t\text{
\ \ \ \ \ \ \ \ \ \ \ \ \ \ \ \ \ \ \ \ \ }(1)\\
-\frac{2H}{W^{2}}=\lambda_{2}a(t)\text{ \ \ \ \ \ \ \ \ \ \ \ \ \ \ \ \ }(2)\\
-\frac{\left(  t+aa^{\prime}\right)  H}{W^{2}}=\lambda_{3}%
s\text{\ \ \ \ \ \ \ \ \ \ \ \ \ \ \ }(3)
\end{array}
\right.  $

Therefore, the problem of classifying the ruled surfaces $r(t,s)=(t,a(t),s)$%
\ by a straight geodesics in $\mathbb{H}_{3}$ satisfying $\left(  1.7\right)
$ is reduced to the integration of the ordinary differential equations of
$\left(  3.2\right)  .$

Next we study it according to the constants $\lambda_{i},i=1,2,3.$

The unique case is\textbf{ }$\lambda_{3}=0$ since the two parts of $\left(
3\right)  $\ are independants. The system $(3.2)$ turn to

$(3.2.a)$\qquad$\qquad\qquad\left\{
\begin{array}
[c]{c}%
\frac{2a^{\prime}H}{W^{2}}=\lambda_{1}t\text{
\ \ \ \ \ \ \ \ \ \ \ \ \ \ \ \ \ \ \ }(1)\\
-\frac{2H}{W^{2}}=\lambda_{2}a(t)\text{ \ \ \ \ \ \ \ \ \ \ \ \ \ \ }(2)\\
-\frac{\left(  t+aa^{\prime}\right)  H}{W^{2}}%
=0\text{\ \ \ \ \ \ \ \ \ \ \ \ \ \ \ }(3)
\end{array}
\right.  $

In the last form, $\left(  3\right)  $ implies $H\equiv0$ or $aa^{\prime}=-t$

\textbf{i)} If $H\equiv0,$ then we get $\lambda_{1}=\lambda_{2}=0$ and the
surface $\mathcal{S}_{1}:r(t,s)=(t,a(t),s)$ is minimal.

\textbf{ii)} If $H\neq0$ and $aa^{\prime}+t=0,$ from which we have

\begin{center}
$aa^{\prime}=-t\Longrightarrow a^{2}=-t^{2}+c,a(t)=\pm\sqrt{c-t^{2}}%
,c\in\mathbb{R}^{+\ast},t\in\left[  -c,c\right]  .$
\end{center}

By multiplying $(2)$ by $a^{\prime}$ and using $aa^{\prime}=-t,$ $(3.2.a)$ became

$(3.2.b)$\qquad$\qquad\qquad\left\{
\begin{array}
[c]{c}%
\frac{2a^{\prime}H}{W^{2}}=-\lambda_{1}aa^{\prime}%
\text{\ \ \ \ \ \ \ \ \ \ \ \ \ \ \ \ \ }(1)\\
-\frac{2a^{\prime}H}{W^{2}}=\lambda_{2}aa^{\prime}%
\text{\ \ \ \ \ \ \ \ \ \ \ \ \ \ \ }(2)\\
-\frac{\left(  t+aa^{\prime}\right)  H}{W^{2}}%
=0\text{\ \ \ \ \ \ \ \ \ \ \ \ \ \ \ \ }(3)
\end{array}
\right.  $

and the sum $(1)+(2)=0$ which implies $\lambda_{1}=\lambda_{2}$.

In fact, since $H=\frac{a^{\prime\prime}}{2W^{2}},$ $aa^{\prime}=-t,$
$W^{2}=1+(a^{\prime})^{2},\ (1)$\ turns out to

\begin{center}
$\left(  \frac{1}{1+(a^{\prime})^{2}}\right)  ^{\prime}=\lambda_{1}\left(
a^{2}\right)  ^{^{\prime}}\Longleftrightarrow\frac{1}{c}\left(  c-t^{2}%
\right)  ^{^{\prime}}=\lambda_{1}\left(  c-t^{2}\right)  ^{^{\prime}}$
\end{center}

We obtain $\lambda_{1}=\frac{1}{c}.$ The same happen with $(2).$\ We conclude
$\lambda_{1}=\lambda_{2}=\frac{1}{c}.$

\subsubsection{Theorem}

\textit{Let }$\mathcal{S}_{1}$\textit{\ be a surface ruled by straight lines
as geodesics in the 3-dimensional }$\mathbb{H}_{3}$ \textit{parametrized by}

\begin{center}
$\mathcal{S}_{1}:r(t,s)=\left(  r_{1}(t,s),r_{2}(t,s),r_{3}(t,s)\right)
=(t,a(t),s).$
\end{center}

\textit{Then, }$r_{i}(t,s)$\textit{\ satisfies the equation }$\Delta
r_{i}=\lambda_{i}r_{i},\lambda_{i}\in\mathbb{R},i=1,2,3,$\textit{\ up to
isometry of }$\mathbb{H}_{3},$\textit{\ if and only if }

$1^{\circ})$\textit{\ }$\mathcal{S}_{1}$\textit{\ has zero mean curvature,
}$\mathcal{S}_{1}:r(t,s)=(t,a(t)=\alpha t+\beta,s);\alpha,\beta\in\mathbb{R},$

$2^{\circ})$\textit{\ The Euclidean circular cylinder }

\begin{center}
$\mathcal{S}_{1}:r(t,s)=(t,a(t)=\pm\sqrt{c-t^{2}},s);c\in\mathbb{R}^{+}%
{}^{^{\ast}},t\in\left[  -c,c\right]  ,s\in\mathbb{R}$.
\end{center}

\bigskip

In $\left[  1\right]  ,\left[  9\right]  ,\left[  10\right]  ,\left[
14\right]  $ the authors proved that, the only well known examples of surfaces
of finite type in Euclidean space $\mathbb{E}^{3}$ are the sphere, the
circular cylinder, and naturally, the minimal surfaces. They proved and gave
the circular cylinder in Euclidean space as the only ruled surface.

We observe in the above theorem, that the circular cylinder$,$ is the only
ruled surface by a geodesic lines which stays of finite type in $\mathbb{H}%
_{3}$ besides minimal surfaces in $\mathbb{H}_{3}$\ which satisfies $\Delta
r_{i}=\lambda_{i}r_{i},\lambda_{i}\in\mathbb{R},i=1,2,3.$

\subsection{The second family of surfaces ruled by geodesic lines satisfying
$\Delta\tilde{r}_{i}=\tilde{\lambda}_{i}\tilde{r}_{i}$}

For the second family of ruled surfaces parametrized as ruled surface

\begin{center}
$\mathcal{S}_{2}:\tilde{r}(x,y)=\left(  \tilde{r}_{1}(x,y)=x,\tilde{r}%
_{2}(x,y)=y,\tilde{r}_{3}(x,y)=z\right)  $
\end{center}

where $\left(  x,y,z\right)  $ are themselves expressed in the parameters
$\left(  t,s\right)  $ as

\begin{center}
$\mathcal{S}_{2}:\left(  x(t,s)=t+su(t),y(t,s)=s,z(t,s)=a(t)+\frac{ts}%
{2}\right)  $

$=\left(  t,0,a\left(  t\right)  \right)  +s\left(  u\left(  t\right)
,1,\frac{t}{2}\right)  .$
\end{center}

The equation $x(t,s)=t+su(t)$ defines implicitly the function $t=t(x,y)$ which
implies $\mathcal{S}_{2},$ is locally as the graph of a differential function
$f.$ From these remarks, about $t(x,y)$ and $f(x,y)$ we have

\begin{center}
$dx=\left(  1+su^{\prime}\right)  dt+uds,dy=ds,dz=f_{x}\left(  \left(
1+su^{\prime}\right)  dt+uds\right)  +f_{y}ds$
\end{center}

where

\begin{center}
$t_{x}=\frac{1}{1+su^{\prime}},t_{y}=-ut_{x};f_{x}=\left(  a^{\prime}+\frac
{y}{2}\right)  t_{x},f_{y}=t_{y}\left(  a^{\prime}+\frac{y}{2}\right)
+\frac{t}{2};$
\end{center}

By help with $(2.3)-(2.7),$ the coefficients of the fundamental forms,

\begin{center}
$\left\{
\begin{array}
[c]{c}%
E=1+\left(  \left(  a^{\prime}+\frac{y}{2}\right)  t_{x}+\frac{y}{2}\right)
^{2}\text{ \ \ }\\
G=1+u^{2}\left(  \left(  a^{\prime}+\frac{y}{2}\right)  t_{x}+\frac{y}%
{2}\right)  ^{2}\\
F=-u\left(  \left(  a^{\prime}+\frac{y}{2}\right)  t_{x}+\frac{y}{2}\right)
^{2}.\text{ \ \ }%
\end{array}
\right.  $
\end{center}

The corresponding Laplacian operator for $\mathcal{S}_{2}$ is, for each
component $\tilde{r}_{i}(x,y)$

$(3.3)$\qquad$\Delta\tilde{r}_{i}=-\frac{1}{\sqrt{EG-F^{2}}}\left(  \left(
\frac{G\tilde{r}_{i,x}-F\tilde{r}_{i,y}}{\sqrt{EG-F^{2}}}\right)  _{x}+\left(
\frac{-F\tilde{r}_{i,x}+E\tilde{r}_{i,y}}{\sqrt{EG-F^{2}}}\right)
_{y}\right)  =\tilde{\lambda}_{i}\tilde{r}_{i}$

As in the previous paragraph, for each coordinate $\tilde{r}_{i},i=1,2,3,$ we have

$\bullet$\textit{\ }For $\tilde{r}_{1}(x,y)=x(t,s)=t+yu(t),$

\begin{center}
$\tilde{r}_{1,x}=t_{x}+yu^{\prime}t_{x}=1;\tilde{r}_{1,y}=t_{y}+u+yu^{\prime
}t_{y}=0$

$G\tilde{r}_{1,x}-F\tilde{r}_{1,y}=1+Q^{2},-F\tilde{r}_{1,x}+E\tilde{r}%
_{1,y}=-PQ,$
\end{center}

$(3.3),$ turns out to, for $i=1$

\begin{center}
$\Delta x=-\frac{1}{W}\left(  \left(  \frac{\left(  1+Q^{2}\right)  }%
{W}\right)  _{x}+\left(  \frac{-PQ}{W}\right)  _{y}\right)  $
\end{center}

thus, by help with $Q=-uP$ we have

\begin{center}
$\Delta x=\frac{1}{W^{4}}\left(
\begin{array}
[c]{c}%
PP_{x}+PQ_{y}+P^{3}Q_{y}+P^{3}u^{2}P_{x}-P^{3}u^{3}P_{y}\\
+P^{3}u^{3}Q_{x}-PuP_{y}+PuQ_{x}+2P^{3}uQ_{x}%
\end{array}
\right)  .$
\end{center}

$\bullet$\textit{\ }For $\tilde{r}_{2}(x,y)=y(t,s)=s,$ $\tilde{r}%
_{2,x}=0;\tilde{r}_{2,y}=1,$

\begin{center}
$G\tilde{r}_{2,x}-F\tilde{r}_{2,y}=-PQ,-F\tilde{r}_{2,x}+E\tilde{r}%
_{2,y}=1+P^{2},$
\end{center}

$(3.3)$ turns out for $i=2$

\begin{center}
$\Delta y=-\frac{1}{W}\left(  \left(  \frac{-PQ}{W}\right)  _{x}+\left(
\frac{1+P^{2}}{W}\right)  _{y}\right)  $
\end{center}

and we obtain, with the help of $Q=-uP$

\begin{center}
$\Delta y(t,s)=-\frac{1}{W^{4}}\left(
\begin{array}
[c]{c}%
PP_{y}-PQ_{x}+P^{3}P_{y}-P^{3}Q_{x}+2P^{3}u^{2}P_{y}\\
+P^{3}u^{3}P_{x}+PuP_{x}+PuQ_{y}+P^{3}uQ_{y}%
\end{array}
\right)  .$
\end{center}

$\bullet$\textit{\ }For $\tilde{r}_{3}(x,y)=f(x,y)=a(t)+\frac{ty}{2},$
$\tilde{r}_{3,x}=P-\frac{y}{2},\tilde{r}_{3,y}=Q+\frac{x}{2},$

\begin{center}
$\left\{
\begin{array}
[c]{c}%
G\tilde{r}_{3,x}-F\tilde{r}_{3,y}=P-\frac{1}{2}y-\frac{1}{2}xPQ-\frac{1}%
{2}yQ^{2},\\
-F\tilde{r}_{3,x}+E\tilde{r}_{3,y}=Q+\frac{1}{2}x+\frac{1}{2}yPQ+\frac{1}%
{2}xP^{2}%
\end{array}
\right.  .$
\end{center}

$(3.3)$ turns out for $i=3$

\begin{center}
$\Delta f(x,y)=-\frac{1}{W}\left(  \left(  \frac{P-\frac{1}{2}\left(
y+xPQ+yQ^{2}\right)  }{W}\right)  _{x}+\left(  \frac{Q+\frac{1}{2}\left(
x+yPQ+xP^{2}\right)  }{W}\right)  _{y}\right)  $
\end{center}

A straightforward computation using $\left(  2.4\right)  -\left(  2.8\right)
$ gives

\begin{center}
$=-\frac{1}{2W^{4}}\left(
\begin{array}
[c]{c}%
2P_{x}+2Q_{y}+2P^{2}Q_{y}+2P^{2}u^{2}P_{x}+PxP_{y}-PxQ_{x}+PyP_{x}\\
+PyQ_{y}-P^{3}uy+2P^{2}uP_{y}+2P^{2}uQ_{x}+P^{3}xP_{y}-P^{3}xQ_{x}\\
+P^{3}yQ_{y}-P^{3}u^{3}y-Puy+2P^{3}u^{2}xP_{y}+P^{3}u^{3}xP_{x}+\\
P^{3}u^{2}yP_{x}+PuxP_{x}+PuxQ_{y}+P^{3}uxQ_{y}+P^{3}uyP_{y}+P^{3}uyQ_{x}%
\end{array}
\right)  $
\end{center}

wich give the third equation

\begin{center}
$\Delta f(x,y)=-\frac{1}{2W^{4}}\left(  \left(  P_{x}+P_{y}\right)  \left(
2+Pux+yP\right)  +PW^{2}\left(  x-yu\right)  \right)  $
\end{center}

and the system $\Delta\tilde{r}_{i}=\tilde{\lambda}_{i}\tilde{r}_{i},i=1,2,3$
turns out to

$(3.5)$ $\ \left\{
\begin{array}
[c]{c}%
-\frac{P}{W^{4}}\left(  uW^{2}-\left(  P_{x}+P_{y}\right)  \right)
=\tilde{\lambda}_{1}x\text{
\ \ \ \ \ \ \ \ \ \ \ \ \ \ \ \ \ \ \ \ \ \ \ \ \ \ \ \ \ \ \ \ \ \ \ \ \ \ \ \ \ \ \ }%
{\small (}\tilde{1}{\small )}\\
-\frac{P}{W^{4}}\left(  W^{2}+u(P_{x}+P_{y})\right)  =\tilde{\lambda}%
_{2}%
y\text{\ \ \ \ \ \ \ \ \ \ \ \ \ \ \ \ \ \ \ \ \ \ \ \ \ \ \ \ \ \ \ \ \ \ \ \ \ \ \ \ \ \ \ \ }%
{\small (}\tilde{2}{\small )}\\
-\frac{1}{2W^{4}}\left(  \left(  P_{x}+P_{y}\right)  \left(  2+Pux+yP\right)
+P\left(  x-yu\right)  W^{2}\right)  =\tilde{\lambda}_{3}\tilde{r}%
_{3}(x,y)\text{\ }{\small (}\tilde{3}{\small )}%
\end{array}
\right.  $

or, by help of $\left(  2.7\right)  ,\left(  2.8\right)  $ in the form

$(3.6.a)$ $\ \ \left\{
\begin{array}
[c]{c}%
-\frac{1}{W^{2}}\left(  Pu-2HPW\right)  =\tilde{\lambda}_{1}x\text{\ \ \ \ \ }%
\ \text{\ \ \ \ \ \ \ \ \ \ \ \ \ \ \ \ \ \ \ \ \ \ \ \ \ \ \ \ \ \ \ \ \ \ }%
{\small (}\tilde{1}{\small )}\\
-\frac{1}{W^{2}}\left(  P+2uHPW\right)  =\tilde{\lambda}_{2}%
y\text{\ \ \ \ \ \ \ \ \ \ \ \ \ \ \ \ \ \ \ \ \ \ \ \ \ \ \ \ \ \ \ \ \ \ \ \ \ \ \ \ }%
{\small (}\tilde{2}{\small )}\\
-\frac{1}{W^{2}}\left(  \frac{1}{2}P\left(  x-uy\right)  +HW\left(
2+Py+Pux\right)  \right)  =\tilde{\lambda}_{3}f(x,y)\text{ \ \ \ }%
{\small (}\tilde{3}{\small ).}%
\end{array}
\right.  $

Multiplying by $y$ and by $x,$ ${\small (}\tilde{1}{\small ),}$ and
${\small (}\tilde{2}{\small )}$ respectively in $\left(  3.6.a\right)  $, we have

\begin{center}
$\left\{
\begin{array}
[c]{c}%
-\frac{1}{W^{2}}\left(  yPu-2PyHW\right)  =\tilde{\lambda}_{1}%
xy\ \ \ \ \ \ \ \ \ \ \ \ \ {\small (}\tilde{1}{\small )}^{\prime}\\
-\frac{1}{W^{2}}\left(  Px+2PuxHW\right)  =\tilde{\lambda}_{2}%
yx\ \ \ \ \ \ \ \ \ \ \ \ \ {\small (}\tilde{2}{\small )}^{\prime}%
\end{array}
\right.  $
\end{center}

The left part of ${\small (}\tilde{3}{\small )},$ by help of ${\small (}%
\tilde{1}{\small )}^{\prime}$ and ${\small (}\tilde{2}{\small )}^{\prime}$
turns out to

\begin{center}
$-\frac{1}{W^{2}}\left(  \frac{1}{2}\left(  Px-Puy\right)  +HW\left(
2+Py+Pux\right)  \right)  =$

$-\frac{1}{2W^{2}}\left(  Px+2PuxHW\right)  +\frac{1}{2W^{2}}\left(
Puy-2PyHW\right)  -\frac{2}{W^{2}}HW=\tilde{\lambda}_{3}f(x,y).$
\end{center}

Thus ${\small (}\tilde{3}{\small )},$ became

\begin{center}
$\frac{2H}{W}=\left(  \frac{\tilde{\lambda}_{2}}{2}-\frac{\tilde{\lambda}_{1}%
}{2}\right)  xy-\tilde{\lambda}_{3}f(x,y){\small .}$
\end{center}

Again, the system $\left(  3.6.a\right)  $ becomes

$(3.6.b)$\qquad$\qquad\qquad\left\{
\begin{array}
[c]{c}%
-\frac{P}{W^{2}}\left(  u-2HW\right)  =\tilde{\lambda}_{1}%
x\text{\ \ \ \ \ \ \ \ \ \ \ \ \ \ \ \ \ \ \ \ \ \ \ }{\small (}\tilde
{1}{\small )}\\
-\frac{P}{W^{2}}\left(  1+2uHW\right)  =\tilde{\lambda}_{2}%
y\text{\ \ \ \ \ \ \ \ \ \ \ \ \ \ \ \ \ \ \ \ \ }{\small (}\tilde
{2}{\small )}\\
\frac{4H}{W}=\left(  \tilde{\lambda}_{2}-\tilde{\lambda}_{1}\right)
xy-2\tilde{\lambda}_{3}f(x,y)\ \ \ \ \ \ \ \ \ \ \ \ \ {\small (}\tilde
{3}{\small ).}%
\end{array}
\right.  $

Therefore, as the same in previous paragraph, the problem of classifying the
ruled surfaces by a straight geodesics in $\mathbb{H}_{3}$ satisfying $\left(
1.7\right)  $ is the integration of the ordinary differential equations of
$\left(  3.6.b\right)  .$

Next we study it according to the constants $\tilde{\lambda}_{i},i=1,2,3.$

With the help of $\left(  2.8\right)  ,$ the system $\left(  3.6\right)
$\ turns to

$(3.6.c)$\qquad$\qquad\qquad\left\{
\begin{array}
[c]{c}%
-\frac{P}{W^{3}}\left(  uW-H_{1}\right)  =\tilde{\lambda}_{1}%
x\text{\ \ \ \ \ \ \ \ \ \ \ \ \ \ \ \ \ \ }{\small (}\tilde{1}{\small )}\\
-\frac{P}{W^{3}}\left(  W+uH_{1}\right)  =\tilde{\lambda}_{2}%
y\text{\ \ \ \ \ \ \ \ \ \ \ \ \ \ \ \ \ \ }{\small (}\tilde{2}{\small )}\\
\frac{2}{W^{3}}H_{1}=\left(  \tilde{\lambda}_{2}-\tilde{\lambda}_{1}\right)
xy-2\tilde{\lambda}_{3}f\ \ \ \ \ \ \ \ \ \ \ {\small (}\tilde{3}{\small ).}%
\end{array}
\right.  $

\textbf{Case} \textbf{1}$.$ Let $\tilde{\lambda}_{3}=0.$

\textbf{i)} If $\tilde{\lambda}_{1}=\tilde{\lambda}_{2}.$ ${\small (}\tilde
{3}{\small )}$ implies $H_{1}=0$ say $H\equiv0$. Thus $\mathcal{S}$ is minimal
and consequently $\tilde{\lambda}_{1}=\tilde{\lambda}_{2}=0$.

\textbf{ii)} If\ $\tilde{\lambda}_{1}=0,$ $\tilde{\lambda}_{2}\neq0.$\ From
$(\tilde{1}{\small )}$ we have $P=0$ or $H_{1}=uW$.

$a.1)$ For $P=0$, $(\tilde{2})$ give $\tilde{\lambda}_{2}=0$. The same begins
if $\tilde{\lambda}_{2}=0$ and $\tilde{\lambda}_{1}\neq0.$ We get contradiction.

Thus, we conclude the same as obtained in \textbf{i)}, $\tilde{\lambda}%
_{1}=\tilde{\lambda}_{2}=\tilde{\lambda}_{3}=0.$

$a.2)$ For $H_{1}=uW$ with $P\neq0,$ the equations\ $(\tilde{2})$ and
${\small (}\tilde{3}{\small )}$ in $(3.6.c)$ turns out to

\begin{center}
$\left\{
\begin{array}
[c]{c}%
-\frac{P}{W^{2}}\left(  1+u^{2}\right)  =\tilde{\lambda}_{2}%
y\text{\ \ \ \ \ \ }{\small (}\tilde{2}{\small )}\\
\frac{2u}{W^{2}}=\tilde{\lambda}_{2}%
xy\ \ \ \ \ \ \ \ \ \ \ \ \ \ \ \ \ \ {\small (}\tilde{3}{\small ).}%
\end{array}
\right.  $
\end{center}

The corresponding surface for $\tilde{\lambda}_{2}\neq0,\tilde{\lambda}%
_{1}=\tilde{\lambda}_{3}=0$, is $f\left(  x,y\right)  $ solution of

$\left(  3.7\right)  \qquad\qquad\qquad\qquad\qquad f_{x}=-\frac{2ux}{1+u^{2}%
}-\frac{y}{2}.$

We obtain the same\ as $a.1),$\ if $\tilde{\lambda}_{2}=0,$ $\tilde{\lambda
}_{1}\neq0$ when $P=0.$

$a.3)$ For $uH_{1}=-W$ with $P\neq0,$ the equations\ $(\tilde{1})$ and
${\small (}\tilde{3}{\small )}$ in $(3.6.c)$ turns out to

\begin{center}
$\left\{
\begin{array}
[c]{c}%
-\frac{P}{uW^{2}}\left(  1+u^{2}\right)  =\tilde{\lambda}_{1}%
x\text{\ \ \ \ \ \ \ \ }{\small (}\tilde{1}{\small )}\\
-\frac{2}{uW^{2}}=-\tilde{\lambda}_{1}%
xy\ \ \ \ \ \ \ \ \ \ \ \ \ \ \ \ {\small (}\tilde{3}{\small ).}%
\end{array}
\right.  $
\end{center}

The corresponding surface for $\tilde{\lambda}_{1}\neq0,\tilde{\lambda}%
_{2}=\tilde{\lambda}_{3}=0$, is $f\left(  x,y\right)  $ solution of

$\left(  3.8\right)  \qquad\qquad\qquad\qquad\qquad f_{x}=\frac{2}{y\left(
1+u^{2}\right)  }-\frac{y}{2}.$

\textbf{iii)} If $\tilde{\lambda}_{1}\neq\tilde{\lambda}_{2},\tilde{\lambda
}_{1}\tilde{\lambda}_{2}\neq0.$ We put $\tilde{\lambda}_{2}-\tilde{\lambda
}_{1}=\tilde{\lambda}.$ The equation ${\small (}\tilde{3}{\small )}$ turn to

$\left(  3.9\right)  \qquad\qquad\qquad\qquad f_{xx}+f_{yy}=\frac{1}{2}%
\tilde{\lambda}_{2}xy\left(  1+\left(  f_{x}+\frac{y}{2}\right)  ^{2}\left(
1+u^{2}\right)  \right)  ^{2}.$

\textbf{Case} \textbf{2}$.$ Let $\tilde{\lambda}_{3}\neq0.$ We replay the
system $(3.6.c)$

\textbf{iii)} $a)$ If $\tilde{\lambda}_{1}=0$ we have $P=0,$ then
$\tilde{\lambda}_{2}=0.$\ The same begins if $\tilde{\lambda}_{2}=0,$ we
obtain $\tilde{\lambda}_{1}=0.$ Thus, we conclude $\tilde{\lambda}_{3}%
\neq0,\tilde{\lambda}_{1}=\tilde{\lambda}_{2}=0.$ Since $P=0$ implies
$Q=-uP=0$ and ${\small (}\tilde{3}{\small )}$\ becomes

\begin{center}
$\frac{1}{W^{4}}\left(  f_{xx}+f_{yy}\right)  =0\Longrightarrow-\tilde
{\lambda}_{3}f=0.$
\end{center}

Thus $\tilde{\lambda}_{3}=0.$ Contradiction, see Case 1, i)\textbf{.}

$b)$ Also, if $\tilde{\lambda}_{1}=0$, with $P\neq0,$ from ${\small (}%
\tilde{1}{\small )}$\ we have $H_{1}=uW$ which implies that ${\small (\tilde
{2})}$ became

\begin{center}
$-\frac{P}{W^{2}}\left(  u^{2}+1\right)  =\tilde{\lambda}_{2}y$
\end{center}

then $\tilde{\lambda}_{2}\neq0.$ The equation ${\small (}\tilde{3}{\small )}%
$\ turn out to

$\left(  3.10\right)  \qquad\ \ \ \ \qquad f_{xx}+f_{yy}=\frac{1}{2}\left(
\tilde{\lambda}_{2}xy-2\tilde{\lambda}_{3}f\right)  \left(  1+\left(
f_{x}+\frac{y}{2}\right)  ^{2}\left(  1+u^{2}\right)  \right)  ^{2}.$

$c)$ If $\tilde{\lambda}_{2}=0$, with $P\neq0,$ from ${\small (\tilde{2})}%
$\ we have $uH_{1}=-W$ which implies that ${\small (}\tilde{1}{\small )}$ became

\begin{center}
$-\frac{P}{uW^{2}}\left(  u^{2}+1\right)  =\tilde{\lambda}_{1}x$
\end{center}

then $\tilde{\lambda}_{1}\neq0.$ The equation ${\small (}\tilde{3}{\small )}
$\ turn to

$\left(  3.11\right)  \qquad\qquad\qquad f_{xx}+f_{yy}=-\frac{1}{2}\left(
\tilde{\lambda}_{1}xy+2\tilde{\lambda}_{3}f\right)  \left(  1+\left(
f_{x}+\frac{y}{2}\right)  ^{2}\left(  1+u^{2}\right)  \right)  ^{2}.$

\textbf{iv)} $a)$ If $\tilde{\lambda}_{2}=\tilde{\lambda}_{1}\neq0.\ $The
equation ${\small (}\tilde{3}{\small )}$\ turn to$\ $

$\left(  3.12\right)  \qquad\qquad\qquad f_{xx}+f_{yy}=-\tilde{\lambda}%
_{3}f\left(  1+\left(  f_{x}+\frac{y}{2}\right)  ^{2}\left(  1+u^{2}\right)
\right)  ^{2}.$

\textbf{v)} $a)$ If $\tilde{\lambda}_{2}=\tilde{\lambda}_{1}=\tilde{\lambda
}_{3}=\tilde{\lambda}\neq0.\ $The equation ${\small (}\tilde{3}{\small )}%
$\ turn to$\ $

$\left(  3.13\right)  \qquad\qquad\qquad f_{xx}+f_{yy}=-\tilde{\lambda
}f\left(  1+\left(  f_{x}+\frac{y}{2}\right)  ^{2}\left(  1+u^{2}\right)
\right)  ^{2}$

To solve and to give explicitly the solutions of $\left(  3.7\right)  -\left(
3.13\right)  $\ is, in my opinion, difficult. The solutions were in the forms
of elliptic integrals which are expressed in Legendre forms, see $\left[
16\right]  $ pp. 276-280.

\subsubsection{Theorem}

\textit{Let }$\mathcal{S}_{2}$\textit{\ be a surface ruled by straight lines
as geodesics in the 3-dimensional }$\mathbb{H}_{3}$ \textit{parametrized by}

\begin{center}
$\mathcal{S}_{2}:\tilde{r}(x,y)=\left(  \tilde{r}_{1}(x,y)=x,\tilde{r}%
_{2}(x,y)=y,\tilde{r}_{3}(x,y)=f\left(  x,y\right)  \right)  $
\end{center}

\textit{where }$\left(  x,y,z\right)  $\textit{ are themeselves expressed in
the parameters }$\left(  t,s\right)  $\textit{ as}

\begin{center}
$\mathcal{S}_{2}:\left(  x(t,s)=t+su(t),y(t,s)=s,f\left(
x(t,s),y(t,s)\right)  =a(t)+\frac{ts}{2}\right)  $

$=\left(  t,0,a\left(  t\right)  \right)  +s\left(  u\left(  t\right)
,1,\frac{t}{2}\right)  .$
\end{center}

\textit{Then, }$\tilde{r}_{i}$\textit{\ satisfies the equation }$\Delta
\tilde{r}_{i}=\tilde{\lambda}_{i}\tilde{r}_{i},\tilde{\lambda}_{i}%
\in\mathbb{R},i=1,2,3,$\textit{\ up to isometry of }$\mathbb{H}_{3}%
,$\textit{\ if and only if }

$1^{\circ})$\textit{\ }$\mathcal{S}_{2}$\textit{\ has zero mean curvature,
}$\mathcal{S}_{2}:\tilde{r}(x\left(  t,s\right)  ,y\left(  t,s\right)
)=\left(  t+su(t),s,a(t)+\frac{ts}{2}\right)  $

$2^{\circ})$\textit{\ }$\mathcal{S}_{2}$\ \textit{is the graphs, solutions of
the abvious partial differential equations }$\left(  3.7\right)  -\left(
3.13\right)  .$\textit{ Among those, there are of them, those of the 1-type,
the 2-type and the 3-type surfaces.}

\bigskip

\textbf{Acknowledgement}

I wish to thank all the colleagues members for the LMIA of Mulhouse University
(France) for their cordial reception, for their kind hospitality and
encouragement during 2015/16, my stay at UHA.

\bigskip

Mohammed BEKKAR

D\'{e}partement de math\'{e}matiques, Facult\'{e} des Sciences,

Universit\'{e} d'Oran I, Alg\'{e}rie.

\textit{or}

Laboratoire de Math\'{e}matique-Informatique et Applications (LMIA),

6, rue des Fr\`{e}res Lumi\`{e}re, Universit\'{e} de Haute Alsace,

68200 Mulhouse, France.

\textit{e-mail address: }bekkar\_99@yahoo.fr

\bigskip



\begin{thebibliography}{99}                                                                                               %
\bibitem {1}C. Baikoussis, B.-Y. Chen, L. Verstraelen, \textit{Ruled surfaces
and tubes with finite type Gauss maps, }Tokyo J. Math. 16, (1993), 341-349.

\bibitem {2}M. Bekkar, \textit{Exemples de surfaces minimales dans l'espace de
Heisenberg. }Rend. Sem. Fac. Sci. Univ. Cagliari 61 (1991), 123-130.

\bibitem {3}M. Bekkar, T. Sari, \textit{Surfaces minimales r\'{e}gl\'{e}es
dans l'espace de Heisenberg }$\mathbb{H}_{3}$; Rend. Sem. Mat. Univ. Pol.
Torino, Vol. 50, 3 (1992).

\bibitem {4}M. Bekkar, \textit{Sur un syst\`{e}me d'\'{e}quations aux
d\'{e}riv\'{e}es partielles dans l'espace de Heisenberg,\ }Rend. Sem. Mat.
Univ. Pol. Torino, 59, n$^{\circ}$. 3, (2003),177--184.

\bibitem {5}M. Bekkar, B. Senoussi, \textit{Translation surfaces in the
3-dimensional space satisfying }$\Delta^{III}r_{i}=\lambda_{i}r_{i};$ J. Geom.
103 (2012), 367-374.

\bibitem {6}B.-Y. Chen, \textit{Surfaces of finite type in Euclidean 3-space,
}Bul. Soc. Math. Belg. Ser. B., 39, (1987), 243-254.

\bibitem {7}B.-Y. Chen, \textit{Finite type Submanifolds and Generalizations,}
University of Roma, Roma (1985).

\bibitem {8}B.-Y. Chen, \textit{Total Mean Curvature and Submanifolds of
finite type, }World Scientific Publisher, Singapore, (1984).

\bibitem {9}B.-Y. Chen; F. Dillen; L. Verstraelen; L. Vrancken, \textit{Ruled
surfaces of finite type, }Bull. Austral. Math. Soc., 42, (1990), 447-453.

\bibitem {10}F. Dillen; J. Pas, L. Vertraelen, \textit{On surfaces of finite
type in Euclidean 3-space,} Kodai Math. J. 13 (1990), n$^{\circ}$. 1, 10-21.

\bibitem {11}F. Dillen; L. Verstraelen; G. Zafindratafa, \textit{A
generalization of the translation surfaces of Scherk, }Diff. Geometry in honor
of Radu Rosca, K. U. L. (1991), 107-109.

\bibitem {12}Do Carmo Manfredo P., \textit{Differential geometry of curves and
surfaces, }Prentice-Hall, Inc., Englewood Cliffs, N.J., 1976. viii+503 pp.

\bibitem {13}C. Figueroa; F. Mercuri; R. Pedrosa, \textit{Invariant surfaces
of the Heisenberg groups, }Ann. Mat. Pura Appl. 177 (1999), 173-194.

\bibitem {14}O. J. Garay, \textit{On a certain class of finite type surfaces
of revolution, }Kodai Math. J. 11, 1 (1988) 25- 31.

\bibitem {15}M. Goze; P. Piu, \textit{Classification des m\'{e}triques
invariantes \`{a} gauche sur le groupe de Heisenberg}, Rend. Circ. Mat.
Palermo, 39-2 (1990), 299-306.

\bibitem {16}S. Gradshtein; M. Ryzhyk, \textit{Tables of integrals, series and
products, }Academic Press, (1980).

\bibitem {17}Th. Hangan, \textit{Sur les distributions totalement
g\'{e}od\'{e}siques du groupe nilpotent riemanien }$\mathbb{H}_{2p+1}$, Rend.
del Sem. Univ. Cagliari (Italie), LV 1 (1985), 17--24.

\bibitem {18}J. Inoguchi; R. L\'{o}pez; M. I. Munteanu, \textit{Minimal
translation surfaces in the Heisenberg group; }Geo. Dedicata, 161, (2011), 221-231.

\bibitem {19}R. L\'{o}pez; M. I. Munteanu, \textit{Minimal translation
surfaces in Sol}$_{3}$, J. Math. Soc. Japan 64 (2012), 985--1003

\bibitem {20}F. Mercuri; S. Montaldo; P. Piu; \textit{A Weierstrass
representation formula for minimal surfaces in Heisenberg space }%
$\mathbb{H}_{3}$\textit{ and }$\mathbb{H}_{2}\times\mathbb{R}$\textit{;} Acta
Math. Sinica. (Engl. Ser.) 22 (2006), n$^{\circ}$. 6, 1603--1612.

\bibitem {21}T. Takahashi, \textit{Minimal immersions of Riemannian
manifolds,} J. Math. Soc. Japan, 18 (1966) 380-385.

\bibitem {22}D. W. Yoon; C. W. Lee; M. K. Karacan, \textit{Some translation
surfaces in the 3-dimensional Heisenberg group,} Bull. Korean Math. Soc. 50
(2013), 1329--1343.
\end{thebibliography}
\end{document}